\baselineskip=16pt

\font\twelvebf=cmbx12

\vskip 15mm
\def\Z_+={{\bf Z}_+}
\def\Z{{\bf Z}}
\def\C{{\bf C}}
\def\N{{\bf N}}
\def\Q{{\bf Q}}
\def\e{\eqno}
\def\l{\ldots}
\def\n{\noindent}

\n
{\twelvebf Highest weight representations of the
quantum algebra
U$_h$(gl$_\infty$)}

\leftskip 50pt
\vskip 32pt
\noindent
T. D. Palev\footnote{$^{a)}$}{E-mail: tpalev@inrne.acad.bg}

\noindent
Institute for Nuclear Research and Nuclear Energy, 1784 Sofia,
Bulgaria
\vskip 12pt

\noindent
N.I. Stoilova\footnote{$^{b)}$}{Permanent address: Institute for 
Nuclear
Research and
Nuclear Energy, 1784 Sofia, Bulgaria; E-mail: stoilova@inrne.acad.bg}

\noindent
International Centre for Theoretical Physics, 34100 Trieste, Italy

\vskip 48pt
\noindent
A class of highest weight irreducible representations of the
quantum algebra $U_h(gl_\infty)$ is constructed. Within each
module a basis is introduced and the transformation relations of
the basis under the action of the Chevalley generators are
explicitly written.

\leftskip 0pt
\vskip 48pt

\vskip 12pt

\n
{\bf I. INTRODUCTION}

\bigskip
The Lie algebra $gl_\infty$ and its completion and central
extension $a_{\infty}\;^{1,2}$  play an important role in several
branches of mathematics and physics. These algebras  are of
interest as examples of Kac-Moody Lie algebras of infinite type
$^{1, 3-5}$. They have applications in the theory of nonlinear
equations$^6$, string theory, two-dimensional statistical models$^7$.
 One of us (T.P.) studies new quantum statistics, based
on the above algebras (see {\it Example 2} in Ref. 8 and the
references therein), leading to local currents$^9$ for what are
called $A-$spinor fields. It is natural to expect that the
deformations of these algebras and their representations may
prove useful too.

The quantum analogues of $gl_\infty$ and $a_{\infty}$ in the
sense of Drinfeld$^{10}$, namely $U_h(gl_\infty)$ and
$U_h(a_\infty)$, were worked out by Levendorskii and Soibelman$^{11}$.
These authors have constructed  a class of highest
weight irreducible representations, writing down explicit
expressions for the transformations of the basis under the action
of the algebra generators.

In the present note we announce results on certain highest weight
irreducible representations (irreps) of $U_h(gl_\infty)$. The
$U_h(gl_\infty)$-modules, we study, are labelled by all possible
complex sequences (see the end of the introduction for the notation)
$$
\{M\} \equiv \{M_i \}_{i\in \Z}\in \C^\infty, \;\;{\rm such \;
that} \;\; M_i-M_j\in \Z_+, \quad \forall i<j.      \eqno(1)
$$
The signatures of the $U_h(gl_\infty)$-modules of Levendorskii
and Soibelman$^{11}$  consist of all those sequences $\{M^{(s)}\}$,
$s\in \Z$, from (1), for which $M_i^{(s)}=1$, if $i< s$ and
$M_i^{(s)}=0$ for $i\ge s$.

Here we use essentially the results of Ref. 12, where a class of
highest weight irreps of $gl_\infty$ was obtained. The
corresponding $gl_\infty$-modules are labelled with the same
signatures (1). We show in fact that each $U_h(gl_\infty)$-module
with a signature $\{M\}$ is an appropriate deformation of
a $gl_\infty$-module with the same signature. 
Our results are based
on a proof of a conjecture we did: the replacement of (almost)
all number multiples in the transformation relations of the basis
$\Gamma (\{M\})$  with formal power series according to
$$
x \rightarrow
[x]={q^x-q^{-x}\over {q-q^{-1}}}\in \C[[h]],\quad q=e^{h/2}\e(2)
$$
turns the $gl_\infty$-module with a signature $\{M\}$ into  
a $U_h(gl_\infty)$-module with the same signature.

\smallskip
\smallskip
It is certainly not surprising such a conjecture holds, because
it takes place for all finite-dimensional irreps of $U_h(gl(n))^{13}$
and of the superalgebras $U_h(gl(n/1)),\;n\in \N^{14}$, if
they are written in a Gel'fand-Zetlin basis, and also for all
essentially typical representations of the superalgebras
$U_h(gl(n/m))^{15}.$
On the other hand counter examples are also
available$^{16}.$   We see our main contribution to be in the
nontrivial proof of the conjecture, not in its formulation.

Throughout the paper we use the notation
(most of them standard): 

$\N$ - all positive integers;

$\Z_+$ - all non-negative integers;

$\Z$ - all integers;

$\Q$ - all rational numbers;

$\C$ - all complex numbers;

$\C[[h]]$ - the ring of all formal power series in $h$ over $\C$;

$\Gamma(\{M\})$ - the $C$-basis of a module with a signature
$\{M\}$;

$q=e^{h/2}$.

\vskip 24pt

\n
{\bf II. THE ALGEBRA U$_h$(gl$_\infty$)}  

\bigskip
Following  Ref. 11 we recall the definition of $U_h(gl_\infty)$. It
is the Hopf algebra, which is topologically free module over
$\C[[h]]$ (complete in $h$-adic topology), with generators
$\{E_i, F_i, H_i  \}_{i\in \Z}$, the 
Chevalley generators, and

\smallskip
\noindent
1. Cartan relations:

$$
\eqalign{
& [H_i, H_j]=0, \cr
& [H_i, E_j]=(\delta _{ij}-\delta _{i,j+1})E_j,  \cr
& [H_i, F_j]=-(\delta _{ij}-\delta _{i,j+1})F_j,  \cr
& [E_i, F_j]=\delta_{ij}{{q^{H_i-H_{i+1}}-q^{-H_i+H_{i+1}}}
  \over{q-q^{-1}}}.\cr
}\e(3)
$$
\noindent
2. $E$-Serre relations:

$$
\eqalign{
& E_iE_j=E_jE_i \quad {\rm if} \;\; |i-j|\neq 1,\cr
& E_i^2E_{i+1}-(q+q^{-1})E_iE_{i+1}E_i+E_{i+1}E_i^2=0, \cr
& E_{i+1}^2E_i-(q+q^{-1})E_{i+1}E_iE_{i+1}+E_iE_{i+1}^2=0. \cr
}\e(4)
$$

\n
3. $F$-Serre relations:

$$
\eqalign{
& F_iF_j=F_jF_i \quad {\rm if} \;\; |i-j|\neq 1,\cr
& F_i^2F_{i+1}-(q+q^{-1})F_iF_{i+1}F_i+F_{i+1}F_i^2=0, \cr
& F_{i+1}^2F_i-(q+q^{-1})F_{i+1}F_iF_{i+1}+F_iF_{i+1}^2=0. \cr
}\e(5)
$$

\bigskip
A counity $\varepsilon $, a comultiplication $\Delta $ and an
antipode $S$ are defined as:
$$
\eqalign{
& \varepsilon (E_i)=\varepsilon (F_i)= \varepsilon (H_i)=0,\cr
& \Delta (H_i)=H_i\otimes 1+1\otimes H_i, \cr
& \Delta (E_i)=E_i\otimes q^{(H_i-H_{i+1})/2} +
  q^{(-H_i+H_{i+1})/2} \otimes E_i,\cr
& \Delta (F_i)=F_i\otimes q^{(H_i-H_{i+1})/2} +
  q^{(-H_i+H_{i+1})/2} \otimes F_i,\cr
& S(H_i)=-H_i, \quad S(E_i)=-qE_i, \quad S(F_i)=-q^{-1}F_i.\cr
}\e(6)
$$
Throughout the tensor products are topological, namely the algebraic
tensor products are replaced with their completion in the $h$-adic
topology.

Note that $\{E_i, F_i, H_i  \}_{i\in \N}$ generate a Hopf subalgebra
$U_h(gl_0(\infty))$ 
of $U_h(gl_\infty)$.

\vskip 24pt

\n
{\bf III. REPRESENTATIONS OF U$_h$(gl$_\infty$) }

\bigskip
Let $W$ be a topologically free $\C[[h]]$-module. We recall that a
$\C[[h]]$-homomorphism $\rho: U_h(gl_\infty) \; \rightarrow \;
End\; W $ is a representation of $U_h(gl_\infty)$ in $W$
(equivalently, $W$ is a $U_h(gl_\infty)$-module) provided $\rho$ is
continuous in the $h$-adic topology.

We proceed to define the $U_h(gl_\infty)$-module $V(\{M\})$ with a
highest weight $\{M\}$ and its (topological) basis.
The basis $\Gamma (\{M\})$, called a central basis ($C$-bases),
is formally the same as the one introduced in Ref. 12 for description 
of
the representations of $gl_\infty$. It consists of all $C-$ patterns

$$
(M)\equiv \left[\matrix
{ .., & M_{1-\theta -k},& \ldots, &M_{-1},&M_0, &M_1, &\ldots
&M_{k+\theta -1},...\cr
.., & \ldots & \ldots & \ldots & \ldots &\ldots  &\ldots &\ldots \cr
&M_{1-\theta -k,2k+\theta -1}, & \ldots, &M_{-1,2k+\theta -1},
&M_{0,2k+\theta -1}, & M_{1,2k+\theta -1}, & \ldots  & M_{k+\theta
-1,2k+\theta -1} \cr
&\ldots &\ldots & \ldots & \ldots  &\ldots &\ldots \cr
& & & M_{-1,3}, & M_{03}, &M_{13} \cr
& & & M_{-1,2}, &M_{02} \cr
& & & & M_{01} \cr
}\right], \eqno(7)
$$
where $k\in \N,\; \theta=0,1$.
Each such pattern is an ordered collection
of complex numbers
$$
M_{i, 2k+\theta -1}, \quad \forall k \in \N, \quad \theta =0,1,
\quad
\;\; i \in [-\theta -k+1, k-1]
 \equiv \{-\theta -k+1,-\theta -k+2,\l, k-1 \}, \eqno(8)
$$
which satisfy the conditions:
\smallskip
\noindent
($i$) there exist $N((M)) \in \N$ 
such that
$$
M_{i,2k+\theta -1}=M_i, \;\; \forall k >N((M)), \quad
 \theta =0,1, \;\; i \in [-\theta -k+1, k-1]; \eqno(9)
$$
(ii)
$$
M_{i+\theta -1,2k+\theta }-M_{i, 2k+\theta -1} \in \Z_+ ,
\;
M_{i,2k+\theta -1}-M_{i+\theta , 2k+\theta} \in \Z_+ ,
\;  \forall k \in \N,  \;
 \theta =0,1, \; i \in [1-\theta -k, k-1]. \eqno(10)
$$

\bigskip

Denote by $V_0(\{M\})$ the free $\C[[h]]$-module with generators
$\Gamma (\{M\})$ and let $V(\{M\})$ be its completion in the $h$-adic
topology. $V(\{M\})$ is topologically free $\C[[h]]$-module with
(topological) basis $\Gamma (\{M\})$. In order to give somewhat
more explicit description of $V(\{M\})$ set $V$ to be the (formal)
complex linear envelope of $\Gamma (\{M\})$. Then 
$V_0(\{M\})=V\otimes_\C \C[[h]]$ and $V(\{M\})$ consists of all
formal power series in $h$ with coefficients in $V$:
$$
v=v_0+v_1h+v_2h^2+\ldots \quad (v_0, v_1, v_2, \ldots \in V).   \e(11)
$$
$End\;V(\{M\})$ is a $\C[[h]]$-module, consisting of all $\C[[h]]-
$linear
maps of $V(\{M\})$. If $a\in End\;V(\{M\})$, then
$$
av=av_0+(av_1)h+(av_2)h^2+\ldots. \e(12)
$$
Therefore the transformation of $V(\{M\})$ under the action of
$a\in End\;V(\{M\})$ is completely defined, if  $a$ is defined on
$\Gamma (\{M\})$.

We pass to turn $V(\{M\})$ into 
a  $U_h(gl_{\infty })$  
module.
To this end introduce first some appropriate notation$^{12}.$ 
Denote by $(M)_{\pm \{ j,p\} }$ and $(M)_{\pm \{ j,p\} }^{\pm \{
l,q\} }$ the 
patterns  obtained from the $C$-pattern $(M)$
(7) after the replacements
$$
M_{lq}\rightarrow M_{lq}\pm 1, \quad
M_{jp}\rightarrow M_{jp}\pm 1, \eqno(13)
$$
correspondingly, 
and let
$$
S(j,l;\nu )=\cases {(-1)^\nu & for $j=l$\cr \hskip 0.3cm 1 & for
$j<l$\cr
-1 &for $j>l$\cr},\quad
\theta(i)=\cases {1 & for
$i\geq 0$ \cr  0 & for $i<0$\cr}, \quad
L_{ij}=M_{ij}-i. \eqno(14)
$$
Set moreover
$$
E^0_i=F_i, \quad E^1_i=E_i, \quad i\in \Z. \e(15)
$$

Let $\{\rho(E_i), \rho(F_i), \rho(H_i)\}_{i\in \Z}$ be a
collection of $\C[[h]]$-endomorphisms of $V(\{M\})$, defined on any
$C$-pattern $(M)\in \Gamma (\{M\})$, as follows:
$$
\eqalignno{
& \rho(E_{-1}^{1-\mu})(M)=( [L_{-1,2}-L_{0,1}-\mu]
[L_{0,1}-L_{0,2}+\mu ])^{1/2}(M)_{-(-1)^\mu \{0,1\} }, \;\quad
\mu =0,1\;, & (16)\cr
& &\cr
& \rho(E^{\;\mu}_{(-1)^\nu i-1})(M)=-\sum_{j=1-i-\nu }^{i-1}
\sum_{l=-i}^{i+\nu-1} S(j,l;\nu )& \cr
& &\cr
& \times \Biggl(-{\prod_{k\not= l=-i}^{i+\nu -1}[L_{k,2i+\nu }-
L_{j,2i+\nu -1}
-(-1)^\nu \mu]\prod_{k=1-i}^{i+\nu-2}[L_{k,2i+\nu-2}-L_{j,2i+\nu
-1}-(-1)^\nu
\mu]\over
{\prod_{k\not= j=1-i-\nu}^{i-1}[L_{k,2i+\nu -1}-L_{j,2i+\nu -1}]
[L_{k,2i+\nu -1}-L_{j,2i+\nu-1}+(-1)^{\mu +\nu}]}} & \cr
& &\cr
& \times {\prod_{k=-i-\nu}^{i}[L_{k,2i+\nu +1}-L_{l,2i+\nu }
+(-1)^\nu (1-\mu )]
\prod_{k\not= j=1-i-\nu}^{i-1}[L_{k,2i+\nu-1}-L_{l,2i+\nu }+(-1)^\nu
(1-\mu )]\over
{\prod_{k\not= l=-i}^{i+\nu -1}[L_{k,2i+\nu }-L_{l,2i+\nu }]
[L_{k,2i+\nu }-L_{l,2i+\nu}+(-1)^{\mu +\nu}]}}\Biggr)^{1/2}& \cr
&& \cr
& \times (M)_{-(-1)^{\mu +\nu }\{ j,2i-1+\nu\} }^{-(-1)^{\mu +\nu }\{
l,2i+\nu\} },\quad  i \in {\bf N},\quad  \mu , \nu =0,1\;,&(17) \cr
&& \cr
& \rho(H_{i})(M)=\left(\sum_{j=-|i|}^{|i|+\theta
(i)-1}M_{j,2|i|+\theta(i)} -
\sum_{j=-|i|+1-\theta (i)}^{|i|-1}M_{j,2|i|+\theta
(i)-1}\right) (M), \quad i\in \Z. & (18)\cr
}
$$

\bigskip
\noindent
If a pattern from the right hand side of 
(17) does not
belong to $\Gamma (\{M\})$, i.e., it is not a $C$-pattern, then
the corresponding term has to be deleted. (The coefficients in
front of all such patterns are  undefined, they contain zero
multiples in the denominators. Therefore an equivalent statement
is that all terms with zeros in the denominators have to be
removed). 
With this convention all
coefficients in front of the
$C-$patterns in r.h.s of (16)-(18) 
are well defined as 
elements from $\C[[h]]$.

\bigskip
\noindent
{\bf Proposition 1.} The endomorphisms $\{\rho(E_i), \rho(F_i),
\rho(H_i)  \}_{i\in \Z}$ satisfy Eqs. (3)-(5) with $\rho(E_i),
\rho(F_i), \rho(H_i)$ substituted for $E_i, F_i, H_i $, respectively.

\bigskip
In the nondeformed case ($h \rightarrow 0$) the above
transformation relations define a representation of $gl_\infty$,
if throughout in (16) and (17) the 
"quantum" 
brackets $[\;\;]$ are assumed
to be usual brackets. The nondeformed Eqs. (16)-(18) were derived
in Ref. 12 from the $gl_0(\infty)$ transformation relations of the
Gel'fand-Zetlin basis. The derivation was based on an isomorphism
$\varphi$ of $gl_\infty$ onto $gl_0(\infty)$, defined as
$$
\varphi(E_{ij})=E_{2|i|+\theta(i),2|j|+\theta(j)}, \quad i,j\in
\Z, \e(19)
$$
where $\{E_{ij}\}_{i,j \in \Z}$ and $\{E_{ij}\}_{i,j \in \N}$ are
the Weyl generators of $gl_\infty$ and $gl_0(\infty)$,
respectively.  In the deformed case ($h\ne 0$) $\varphi$ does not
define any more an isomorphism of $U_h(gl_{\infty })$ onto
$U_h(gl_0(\infty))$.  In fact we do not know whether
$U_h(gl_{\infty })$ is isomorphic to its subalgebra
$U_h(gl_0(\infty))$.
The proof of {\it Proposition 1} is based on a
direct verification of the relations (3)-(5). This verification
is extremely lengthy and in certain cases highly nontrivial. The
most difficult 
to check 
are the last Cartan relations in (3),
corresponding to $i=j$. In order to show they hold, one has to
prove as an intermediate step that the following identities are
fulfilled:

\bigskip
$$
\eqalignno{
&  {\scriptstyle {\rm \sum_{j=1-k}^{k-1}\sum_{l=-k}^{k-1}
{\prod_{i\neq l=-k}^{k-1}[L_{i,2k}-L_{j,2k-1}-1]
\prod_{i=1-k}^{k-2}[L_{i,2k-2}-L_{j,2k-1}-1]
\prod_{i=-k}^k[L_{i,2k+1}-L_{l,2k}]
\prod_{i\neq j=1-k}^{k-1}[L_{i,2k-1}-L_{l,2k}]
\over {\prod_{i\neq j=1-k}^{k-1}[L_{i,2k-1}-L_{j,2k-1}]
[L_{i,2k-1}-L_{j,2k-1}-1]
\prod_{i\neq l=-k}^{k-1}[L_{i,2k}-L_{l,2k}]
[L_{i,2k}-L_{l,2k}-1]}}}} & \cr
&&\cr
&&\cr
& {\scriptstyle {\rm -\sum_{j=1-k}^{k-1}\sum_{l=-k}^{k-1}
{\prod_{i\neq l=-k}^{k-1}[L_{i,2k}-L_{j,2k-1}]
\prod_{i=1-k}^{k-2}[L_{i,2k-2}-L_{j,2k-1}]
\prod_{i=-k}^k[L_{i,2k+1}-L_{l,2k}+1]
\prod_{i\neq j=1-k}^{k-1}[L_{i,2k-1}-L_{l,2k}+1]
\over {\prod_{i\neq j=1-k}^{k-1}[L_{i,2k-1}-L_{j,2k-1}]
[L_{i,2k-1}-L_{j,2k-1}+1]
\prod_{i\neq l=-k}^{k-1}[L_{i,2k}-L_{l,2k}]
[L_{i,2k}-L_{l,2k}+1]}}}} & \cr
&&\cr
&&\cr
& {\scriptstyle {\rm =\left[\sum_{j=-k+1}^{k-1}L_{j,2k-1}-
\sum_{j=-k+1}^{k-2}L_{j,2k-2}-
\sum_{j=-k}^k L_{j,2k+1}+\sum_{j=-k}^{k-1} L_{j,2k}-1\right],
\quad k\in \N}}, &
(20)
\cr
}
$$

\bigskip
$$
\eqalignno{
& {\scriptstyle {\rm \sum_{j=-k}^{k-1}\sum_{l=-k}^{k}
{\prod_{i\neq l=-k}^{k}[L_{i,2k+1}-L_{j,2k}+1]
\prod_{i=1-k}^{k-1}[L_{i,2k-1}-L_{j,2k}+1]
\prod_{i=-k-1}^k[L_{i,2k+2}-L_{l,2k+1}]
\prod_{i\neq j=-k}^{k-1}[L_{i,2k}-L_{l,2k+1}]
\over {\prod_{i\neq j=-k}^{k-1}[L_{i,2k}-L_{j,2k}]
[L_{i,2k}-L_{j,2k}+1]
\prod_{i\neq l=-k}^{k}[L_{i,2k+1}-L_{l,2k+1}]
[L_{i,2k+1}-L_{l,2k+1}+1]}}}} & \cr
&&\cr
&&\cr
& {\scriptstyle {\rm -\sum_{j=-k}^{k-1}\sum_{l=-k}^{k}
{\prod_{i\neq l=-k}^{k}[L_{i,2k+1}-L_{j,2k}]
\prod_{i=1-k}^{k-1}[L_{i,2k-1}-L_{j,2k}]
\prod_{i=-k-1}^{k}[L_{i,2k+2}-L_{l,2k+1}-1]
\prod_{i\neq j=-k}^{k-1}[L_{i,2k}-L_{l,2k+1}-1]
\over {\prod_{i\neq j=-k}^{k-1}[L_{i,2k}-L_{j,2k}]
[L_{i,2k}-L_{j,2k}-1]
\prod_{i\neq l=-k}^{k}[L_{i,2k+1}-L_{l,2k+1}]
[L_{i,2k+1}-L_{l,2k+1}-1]}}}} & \cr
&&\cr
&&\cr
& {\scriptstyle {\rm =\left[\sum_{j=-k-1}^{k}L_{j,2k+2}-
\sum_{j=-k}^{k}L_{j,2k+1}-
\sum_{j=-k}^{k-1} L_{j,2k}+\sum_{j=-k+1}^{k-1} L_{j,2k-1}-1\right],
\quad k\in \N }}.
&(21) \cr
}
$$
We prove these identities (and several others, which are
simpler) in a paper to come, where we extend our results to the
completion and central extension $U_h(a_\infty)$ of
$U_h(gl_\infty)$.

A (topological) basis $\{e_n\}_{n\in \N}$ in $U_h(gl_\infty)$ was
given in Ref. 11. Each basis vector $e_n$ is an appropriate
$\C[[h]]$-polynomial in the Chevalley generators. The
$\C[[h]]$-span $W[[h]]$ of the basis is dense in
$U_h(gl_\infty)$. It consists of all $\C[[h]]$-polynomials in the
Chevalley generators.  Extend the domain of definition of $\rho$
on $W[[h]]$ in a natural way: if $\rho$ has already been defined
on $a, b \in W[[h]]$, then set
$$
\rho(\alpha a + \beta b)=\alpha \rho(a) + \beta \rho(b),\quad
\rho(ab)=\rho(a)\rho(b), \quad a, b \in W[[h]], \quad
\alpha, \beta \in \C[[h]]. \e(22)
$$
Denote by $U(gl_\infty)$ the $\C-$linear span of the basis.
Clearly $U(gl_\infty) \subset W[[h]]$ and therefore $\rho$
is defined on $U(gl_\infty)$.
$U_h(gl_\infty)$ consists of all elements of the form
$$
a=\sum_{i=0}^\infty a_i h^i, \quad
(a_0, a_1, a_2, \ldots \in U(gl_\infty)).
\e(23)
$$
For any $i$, according to (11),
$$
\rho(a_i)=\sum_{j=0}^\infty \alpha_{ij}h^j \in V(\{M\}),\quad
\alpha_{ij}\in V.\e(24)
$$
Therefore,
$$
\sum_{i=0}^\infty \rho(a_i)h^i=\sum_{i=0}^\infty
\left(\sum_{j=0}^\infty \alpha_{ij}h^j\right) h^i
=\sum_{n=0}^\infty \left(\sum_{m=0}^n \alpha_{n-m,m}\right)h^n
\in V(\{M\}).\e(25)
$$
Using (25), extend $\rho$ on $U_h(gl_\infty)$:
$$
\rho(a)=\sum_{i=0}^\infty \rho(a_i)h^i\in V(\{M\})\;\;\; \forall \;\;
a\in U_h(gl_\infty).
\e(26)
$$
Hence $\rho$ is a well defined map from $U_h(gl_\infty)$ into
$End\;V(\{M\})$,
$$
\rho: \; U_h(gl_\infty) \rightarrow End\;V(\{M\}). \e(27)
$$

Let $a\in U_h(gl_\infty)$. Then any neighbourhood $W(\rho(a))$
of $\rho(a)$ contains a basic neighbourhood
$\rho(a)+h^nEnd\;V(\{M\})\subset W(\rho(a)) $.
We have in mind the
$h$-adic topology both in $U_h(gl_\infty)$ and $End\;V(\{M\})$.
Evidently
$$
\rho\left( a+h^nU_h(gl_\infty) \right) \subset
\rho(a)+h^nEnd\;V(\{M\})\subset W(\rho(a)) \;\;\forall a\in 
U_h(gl_\infty).
\e(28)
$$
Therefore $\rho$ is a continuous map. 
It satisfies Eq. (22) for any $a, b \in U_h(gl_\infty)$ and
$\alpha, \beta \in \C[[h]]$. Therefore $\rho$ is 
a $\C[[h]]$-homomorphism of $U_h(gl_\infty)$ in $End\;V(\{M\})$.
We have obtained the following
result.

\bigskip
\noindent
{\bf Proposition 2.} The map (27), acting on the $C-$basis according
to Eqs. (16)-(18), defines a representation of $U_h(gl_\infty)$ in
$V(\{M\})$.

\smallskip
Any $U_h(gl_\infty)$-module $V(\{M\})$ is a highest weight module
with respect to the "Borel" subalgebra $N_+$, consisting of all
$\C[[h]]$-polynomials of the unity and $\{E_i\}_{i\in \Z}$.  The
highest weight vector $(\hat{M})$, which by definition satisfies the
condition $\rho(N_+)(\hat{M})=0$, corresponds to the one from (7) with
$$
\hat{M}_{i, 2k+\theta -1}=M_i , \quad \forall k \in \N, \quad \theta 
=0,1,
\quad
\;\; i \in [-\theta -k+1, k-1]. \eqno(29)
$$
Moreover, $V(\{M\})=\rho(U_h(gl_\infty))(\hat{M})$.  Since
$V(\{M\})$ contains no other singular vectors, vectors
annihilated from $\rho(N_+)$, each $V(\{M\})$ is an irreducible
$U_h(gl_\infty)$-module. The proof of the latter follows from the
results in Ref. 12 and the observation that each (deformed)
matrix element in the transformation relations (16)-(18) is zero
only if the corresponding nondeformed matrix element vanishes.

\smallskip
We note in conclusion that all our results remain valid for
$h\notin i\pi\Q$, namely in the case $q$ is 
a number, which is  not a root of 1.

\vskip 48pt
\n
{\bf ACKNOWLEDGMENTS}

\bigskip
One of us (N.I.S.) is grateful to Prof. M.D. Gould for the
invitation to work in his group at the Department of Mathematics
in The University of Queensland,
where part of the work on the present investigation was
completed. 
The stimulating discussions with him are greatly acknowledged.
N.I.S. is thankful to Prof.  Randjbar-Daemi for the kind
hospitality at the High Energy Section of ICTP. This work was
supported by the Australian Research Council and by the Grant
$\Phi -416$ of the Bulgarian Foundation for Scientific Research.

\bigskip
\bigskip

\vskip 12pt
{\settabs \+  $^{11}$ & I. Patera, T. D. Palev, Theoretical
   interpretation of the experiments on the elastic \cr

\+ $^1$ & V.G. Kac  and D.H. Peterson,   Proc. Natl. Acad. Sci.
           USA {\bf 78}, 3308 (1981). \cr

\+ $^2$ & E. Date, M. Jimbo, M. Kashiwara and T. Miwa,  
         J. Phys. Soc. Japan {\bf 50}, 3806 (1981).\cr

\+ $^3$ & V.G. Kac,   {\it Infinite dimensional Lie algebras} {\bf 44}
           (Cambridge: CUP, 1985). \cr

\+ $^4$ & V.G. Kac and A.K. Raina, {\it Bombay lectures on highest
           weight representations of infinite-}  \cr
\+       &  {\it dimensional Lie algebras} in Advanced Series in
           Mathematics {\bf 2} (World Scientific, Singapore, 1987). 
\cr

\+ $^5$ & B. Feigin, D. Fuchs,  {\it Representations of the Virasoro
           algebra} in Representations of infinite- \cr
\+       & dimensional Lie groups and Lie algebras (New York: Gorgon 
and
           Breach, 1989). \cr

\+ $^6$ & E. Date, M. Jimbo, M. Kashiwara and T. Miwa, {\it
           Transformation groups for soliton equations} \cr
\+       & Publ. RIMS Kyoto Univ. {\bf 18}, 1077 (1982).\cr

\+ $^7$ & P. Goddard and D. Olive,  Int. J. Mod. Phys. {\bf A
           1}, 303 (1986).\cr

\+ $^8$ & T.D. Palev,  Rep. Math. Phys. {\bf 31}, 241 (1992).\cr

\+ $^9$ & T.D. Palev,  C.R. Acad. Bulg. Sci. {\bf 32}, 159 (1979).\cr

\+ $^{10}$ & V. Drinfeld,  {\it Quantum groups,} ICM proceedings
            Berkeley 798 (1986).\cr

\+ $^{11}$ & S. Levendorskii and Y. Soibelman,  Comm. Math. Phys.
            {\bf 140}, 399 (1991).\cr

\+ $^{12}$ & T.D. Palev,   J. Math. Phys. {\bf 31}, 579
             (1990); see also
             Funkt. Anal. Prilozh. {\bf 24}, \# 1, 82
             (1990) and \cr
\+        & Funct. Anal. Appl. {\bf 24}, 72 (1990) (English
            transl.). \cr

\+ $^{13}$ & M. Jimbo,   Lect. Notes in Physics,  Berlin,
            Heidelberg, New York: Springer {\bf 246}, 334 (1986).\cr

\+ $^{14}$ & T.D. Palev and V.N. Tolstoy, Comm. Math. Phys.
            {\bf 141}, 549 (1991).\cr

\+ $^{15}$ & T.D. Palev, N.I. Stoilova and J. Van der Jeugt,
         Comm. Math. Phys. {\bf 166}, 367 (1994).\cr

\+ $^{16}$ & Ky Anh Nguyen  and N.I. Stoilova,   J. Math.
          Phys. {\bf 36}, 5979 (1995).\cr

\end